\documentclass[12pt]{amsart}
\usepackage{latexsym}
\usepackage{amsfonts}
\usepackage{isolatin1}
\newtheorem{lemma}{Lemma}
\newtheorem{theorem}{Theorem}
\newcommand{\ovprt}{\overline{\partial}}

\begin{document}
\title{Schrödinger operators with magnetic fields and the canonical solution operator to $\ovprt .$}
\author{Friedrich Haslinger}
\date{}
\address{Institut f\"ur Mathematik, Universit\"at Wien,
Strudlhofgasse 4, A-1090 Wien, Austria;
friedrich.haslinger@univie.ac.at \hskip 0.6cm
\newline http://www.mat.univie.ac.at/\~\,has/ }

\begin{abstract} In this paper we characterize compactness of the canonical solution operator
to $\ovprt $ on weigthed $L^2$ spaces on $\mathbb C.$
For this purpose we consider certain Schrödinger
operators with magnetic fields and use a condition which is equivalent to the
property that these operators have compact resolvents.

\end{abstract}
\maketitle

{\bf 1. Introduction.}

\vskip 1 cm
A nonnegative Borel measure $\nu $ defined on $\mathbb C$ is said to be doubling if
there exists a constant $C$ such that for all $z\in \mathbb C$ and $r\in \mathbb R^+,$
$$\nu (B(z,2r))\le C \nu (B(z,r)).$$
$\mathcal D$ denotes the set of all doubling measures $\nu $ for which there exists
a constant $\delta $ such that for all $z\in \mathbb C,$
$$\nu (B(z,1))\ge \delta.$$
Let $\varphi : \mathbb C\longrightarrow \mathbb R$ be a subharmonic function.
Then $\Delta \varphi $ defines a nonnegative Borel measure, which is finite on compact sets.

Let $\mathcal W$ denote the set of all subharmonic functions $\varphi : \mathbb C
\longrightarrow \mathbb R$ such that $\Delta \varphi \in \mathcal D.$
For $\varphi \in \mathcal W$ we consider the Hilbert spaces 
$$L^2_{\varphi }=\{ f: \mathbb C \longrightarrow \mathbb R \  \mbox{measureable} \ : \ 
\|f\|^2_{\varphi }:=\int_{\mathbb C}|f(z)|^2 \,e^{-2\varphi (z)}\,d\lambda (z)<\infty \}.$$

M. Christ  \cite{Ch} proved that for $\varphi \in \mathcal W$ and for every $f\in L^2_{\varphi }$
there exists $u\in L^2_{\varphi }$ satisfying $\ovprt u=f.$ 
In fact it is shown that there exists a continuous solution operator $\tilde S:L^2_{\varphi}
\longrightarrow L^2_{\varphi }$ for $\ovprt ,$ i.e. $\|\tilde S(f) \|_{\varphi }\le C
\|f\|_{\varphi }$ and $\ovprt \tilde S (f)=f.$
A further result from there indicates that
for $\varphi \in \mathcal W$ the space $A^2_{\varphi }=\mbox{ker}\, \ovprt $ 
of entire functions in $L^2_{\varphi }$
has infinite dimension.

Let $P_{\varphi }:L^2_{\varphi }\longrightarrow A^2_{\varphi }$ denote the Bergman projection.
Then $S=(I-P_{\varphi })\tilde S $ is the uniquely determined canonical solution operator to $\ovprt ,$
i.e. $\ovprt S(f)=f $ and $S(f)\perp A^2_{\varphi }.$

In this paper we characterize compactness of the canonical solution operator to $\ovprt $
on $L^2_{\varphi }$ using 
results about Schrödinger operators with magnetic fields in $\mathbb R^2.$
In an earlier paper we  showed that the canonical solution operator to $\ovprt $
restricted  to  $(0,1)$-forms with holomorphic coefficients can
be  expressed by an integral operator using the Bergman kernel.
This  result can for instance be used to prove that in the case of the unit disc
in $\mathbb C ,$ the canonical solution operator to $\ovprt $
restricted   to  $(0,1)$-forms with holomorphic coefficients is
a Hilbert-Schmidt operator (\cite{Has1} ). With similar methods one can prove that 
the canonical solution operator to $\ovprt $ restricted to the Fock space $A^2_{\varphi },$
where $\varphi (z)=|z|^2$ fails to be compact, whereas in the case $\varphi (z)=|z|^m \ , \ 
m>2$ the canonical solution operator to $\ovprt $ restricted to $A^2_{\varphi }$ is compact
but fails to be Hilbert Schmidt (\cite{Has2}). 

The question of compactness of the solution operator to $\ovprt $ is of interest for various 
reasons - see \cite{FS1} and \cite{FS2} for an excellent survey and \cite{C}, \cite{CD},
\cite{K}, \cite{L}.

A similar situation appears in \cite{SSU} where the Toeplitz $C^*$ -algebra
$\mathcal T (\Omega )$ is considered and the relation between the structure of
$\mathcal T (\Omega )$ and the $\ovprt $-Neumann problem is discussed
(see \cite{SSU} ).

The connection of $\ovprt $ with the theory of Schrödinger operators with magnetic fields
appears in \cite{Ch}, \cite{B} and \cite{FS3}.

Multiple difficulties arise in the case of several complex variables, mainly because the 
geometric structures underlying the analysis become much more complicated.
\vskip 1 cm

{\bf 2. Schrödinger operators with magnetic fields.}

\vskip 1 cm
We want to solve $\ovprt u=f$ for $f\in L^2_{\varphi }.$ The canonical solution
operator to $ \ovprt $ gives a solution with minimal $L^2_{\varphi }$-norm.
Following \cite{Ch} we substitute $v=u\,e^{-\varphi }$ and $g=f\,e^{-\varphi }$
and the equation becomes
$$\overline D v=g \ , \ \mbox{where} \ \overline D = e^{-\varphi }\, \frac {\partial }
{\partial \overline z}\, e^{\varphi }.$$
$u$ is the minimal solution to the $\ovprt $-equation in $L^2_{\varphi }$ if and only if
$v$ is the solution to $\overline D v=g$ which is minimal in $L^2(\mathbb C ).$

The formal adjoint of $\overline D$ is $D=-e^{\varphi }\frac{\partial}{\partial z}e^{
-\varphi }.$ As in \cite{Ch} we define $\mbox{Dom}(\overline D)=\{f\in L^2(\mathbb C)
\ : \ \overline D f\in L^2(\mathbb C ) \}$ and likewise for $D.$ Then $\overline D $ and
$D$ are closed unbounded linear operators from $L^2(\mathbb C)$ to itself.
Further we define $\mbox{Dom}(\overline D D)=\{ u\in \mbox{Dom}(D) \ : \ 
Du\in \mbox{Dom}(\overline D )\}$ and we define $\overline D D $ as $\overline D \circ
D$ on this domain. Any function of the form $e^{\varphi }\,g,$ with $g\in \mathcal C^2_0$
belongs to $\mbox{Dom}(\overline D D)$ and hence $\mbox{Dom}(\overline D D)$ is dense
in $L^2(\mathbb C).$ Since $\overline D=\frac{\partial }{\partial \overline z}+
\frac{\partial \varphi }{\partial \overline z}$ and $D=-\frac{\partial }{\partial z}+
\frac{\partial \varphi }{\partial z}$ we see that
$$\overline D D=-\frac{\partial^2}{\partial z \partial \overline z}-\frac{\partial
\varphi }{\partial \overline z}\,\frac{\partial}{\partial z}+ \frac{\partial \varphi }
{\partial z}\,\frac{\partial}{\partial \overline z}+\left |\frac{\partial \varphi }
{\partial z} \right |^2+\frac{\partial^2 \varphi}{\partial z \partial \overline z}$$
$$=-\frac{1}{4}\, ((d-iA)^2 - \Delta \varphi ),$$
where $A=A_1\,dx+A_2\,dy=-\varphi_y\,dx + \varphi_x\,dy.$ Hence $\overline D D$ is a Schr\"odinger
operator with electric potential $\Delta \varphi $ and with magnetic field $B=dA.$

Now let $\|u\|^2=\int_{\mathbb C}|u(z)|^2\,d\lambda (z)$ for $u\in L^2(\mathbb C)$ and
$(u,v)=\int_{\mathbb C}u(z)\overline{v(z)}\,d\lambda (z)$ denote the inner product of
$L^2(\mathbb C).$

For $\phi , \psi \in \mathcal C^{\infty}_0 (\mathbb C) $ we set
$$h_{A,\varphi}(\phi ,\psi )=(-\frac{1}{4}((d-iA)^2-\Delta \varphi)\phi ,\psi ) =
\sum_{j=1}^2(\frac{1}{4}(\Pi_j(A)\phi, \Pi_j(A)\psi)+((\Delta \varphi)\phi,\psi),$$
where
$$\Pi_j(A)=\frac{1}{i}(\frac{\partial}{\partial x_j}-A_j) \ , \ j=1,2,$$
hence $h_{A,\varphi}$ is a nonnegative symmetric form on $\mathcal C^{\infty}_0(\mathbb C).$

In \cite{Ch} the following results are proved

\begin{lemma} Let $\varphi \in \mathcal W.$
If $u\in \mbox{Dom}(D)$ and $Du\in \mbox{Dom}(\overline D),$ then
$$\|Du\|^2= (\overline D(Du),u).$$
$\overline D D$ is a closed operator and
$$\|u\| \le C\|\overline D D u\|$$
for all $u\in \mbox{Dom}(\overline D D).$ Moreover, for any $f\in L^2(\mathbb C)$
there exists a unique $u\in \mbox{Dom}(\overline D D)$ satisfying $\overline D D u=f.$
Hence $(\overline D D )^{-1}$ is a bounded operator on $L^2(\mathbb C).$
\end{lemma}

The closure $\overline h_{A, \varphi}$ of the form $h_{A,\varphi }$ is a nonnegative symmetric
form. The selfadjoint operator associated with $\overline h_{A,\varphi }$ is the operator
$\overline D D $ from the above Lemma 1 (see \cite{I} and \cite{CFKS}).

\begin{lemma} Let $\varphi \in \mathcal W.$
The canonical solution operator $S: L^2_{\varphi }\longrightarrow
L^2_{\varphi }$ to $\ovprt $ is compact if and only if $(\overline D D)^{-1} : L^2(\mathbb C)
\longrightarrow L^2(\mathbb C )$ is compact.
\end{lemma}

\begin{proof}  If $v$ is the minimal solution to $\ovprt v=g$ in $L^2_{\varphi }$ then
$u=v\,e^{-\varphi }$  is the minimal solution to $\overline D u =g\,e^{-\varphi }$
in $L^2(\mathbb C).$ Hence the canonical solution operator $S$ to $\ovprt $ is compact
if and only if the canonical solution operator to $\overline D u =f $ is compact.
By Lemma 1. we have
$$\|D(\overline D D)^{-1}f\|^2=(\overline D D(\overline D D)^{-1}f, (\overline D D)^{-1}f)=
(f,(\overline D D)^{-1}f)\le \|(\overline D D )^{-1}\|\|f\|^2, $$
hence
$$\|D(\overline D D)^{-1}f\| \le \|(\overline D D )^{-1}\|^{1/2} \, \|f\|$$
and $T=D(\overline D D)^{-1}$ is a bounded operator on $L^2(\mathbb C)$ with
$\overline D Tf=f$ and $Tf\perp \mbox{ker}\overline D ,$ which means that
$T$ is the canonical solution operator to $\overline D u=f.$ Since $(\overline D D)^{-1}$
is a selfadjoint operator (see for instance \cite{I}) it follows that
$$(\overline D D)^{-1}=T^*T.$$
Now we are finished, since $T$ is compact if and only if $T^*T$ is compact (see \cite{W}).
\end{proof}

Using the main theorem in \cite{I} we get

\begin{theorem}Let $\varphi \in \mathcal W.$ The canoical solution operator
$S:L^2_{\varphi}\longrightarrow L^2_{\varphi }$ to $\ovprt $ is compact if and only if
there exists a real valued continuous function $\mu $ on $\mathbb C$ such that
$\mu (z)\to \infty$ as $|z|\to \infty $ and
$$h_{A,\varphi }(\phi, \phi )\ge \int_{\mathbb C}\mu (z)\, |\phi (z)|^2\,d\lambda (z)$$
for all $\phi \in \mathcal C^{\infty}_0(\mathbb C).$
\end{theorem}
\vskip 0.5 cm

\begin{proof} In \cite{I} it is shown that the operator $\overline D D$ has compact resolvent
if and only if the condition in Theorem 1 holds. Therefore, by Lemma 2, Theorem 1 is proved.
\end{proof}

\newpage

\begin{theorem} If $\varphi (z)=|z|^2,$ then the canonical solution operator
$S:L^2_{\varphi}\longrightarrow L^2_{\varphi }$ to $\ovprt $
fails to be compact.
\end{theorem}

\begin{proof} In our case the magnetic field $B$ is the form $B=dA=B(x,y) dx\wedge dy=
\Delta \varphi dx\wedge dy.$
Hence for $\varphi (z)=|z|^2$ we have $\Delta \varphi (z)=4 $ for each $z \in \mathbb C.$
Let $Q_w $ be the ball centered at $w$ with radius $1.$ Then
$$\int_{Q_w}(|B(x,y)|^2+\Delta \varphi (z))\,d\lambda (z) $$
is a constant as $|w|\to \infty ,$ so the assertion follows from \cite{I} Theorem 5.2.

\end{proof}

\begin{theorem} Let $\varphi \in \mathcal W$ and suppose that $\Delta \varphi (z)\to \infty$
as $|z|\to \infty .$ Then the
canonical solution operator $S:L^2_{\varphi}\longrightarrow L^2_{\varphi }$ to $\ovprt $ is
compact.
\end{theorem}

\begin{proof} Since in our case $|B(x,y)|=\Delta \varphi (z)\to \infty,$ as $|z|\to \infty $
the conclusion follows from Lemma 2. and \cite{AHS}, \cite{D} or \cite{I}.
\end{proof}

\vskip 0.5 cm

{\bf Remark.} In \cite{Has2} it shown that for $\varphi (z) =|z|^2$
even the restriction of the canonical solution
$S$ to the Fock space $A^2_{\varphi }$ fails to be compact and that for $\varphi (z)=|z|^m \ ,
\ m>2$ the restriction of $S$ to $A^2_{\varphi }$ fails to be Hilbert Schmidt.

\end{document}